\documentclass[12pt]{article}
\usepackage{amssymb,
amsfonts,epic,graphics,pictex}
\begin{document}

\newtheorem{lem}{Lemma}[section]
\newtheorem{prop}{Proposition}
\newtheorem{con}{Construction}[section]
\newtheorem{defi}{Definition}[section]
\newtheorem{coro}{Corollary}[section]
\newcommand{\hf}{\hat{f}}
\newtheorem{fact}{Fact}[section]
\newtheorem{theo}{Theorem}
\newcommand{\Br}{\Poin}
\newcommand{\Cr}{{\bf Cr}}
\newcommand{\dist}{{\bf dist}}
\newcommand{\diam}{\mbox{diam}\, }
\newcommand{\mod}{{\rm mod}\,}
\newcommand{\compose}{\circ}
\newcommand{\dbar}{\bar{\partial}}
\newcommand{\Def}[1]{{{\em #1}}}
\newcommand{\dx}[1]{\frac{\partial #1}{\partial x}}
\newcommand{\dy}[1]{\frac{\partial #1}{\partial y}}
\newcommand{\Res}[2]{{#1}\raisebox{-.4ex}{$\left|\,_{#2}\right.$}}
\newcommand{\sgn}{{\rm sgn}}
\newcommand{\ph}{{\varphi}}

\newcommand{\C}{{\bf C}}
\newcommand{\D}{{\bf D}}
\newcommand{\Dm}{{\bf D_-}}
\newcommand{\N}{{\bf N}}
\newcommand{\R}{{\bf R}}
\newcommand{\Z}{{\bf Z}}
\newcommand{\tr}{\mbox{Tr}\,}
\newcommand{\h}{{\bf H}}

\newenvironment{nproof}[1]{\trivlist\item[\hskip \labelsep{\bf Proof{#1}.}]}
{\begin{flushright} $\square$\end{flushright}\endtrivlist}
\newenvironment{proof}{\begin{nproof}{}}{\end{nproof}}

\newenvironment{block}[1]{\trivlist\item[\hskip \labelsep{{#1}.}]}{\endtrivlist}
\newenvironment{definition}{\begin{block}{\bf Definition}}{\end{block}}

\newtheorem{conjec}{Conjecture}
\newtheorem{com}{Comment}
\newtheorem{exa}{Example}
\font\mathfonta=msam10 at 11pt
\font\mathfontb=msbm10 at 11pt
\def\Bbb#1{\mbox{\mathfontb #1}}
\def\lesssim{\mbox{\mathfonta.}}
\def\suppset{\mbox{\mathfonta{c}}}
\def\subbset{\mbox{\mathfonta{b}}}
\def\grtsim{\mbox{\mathfonta\&}}
\def\gtrsim{\mbox{\mathfonta\&}}

\newcommand{\ar}{{\bf area}}
\newcommand{\1}{{\bf 1}}
\newcommand{\e}{\varepsilon}
\newcommand{\Int}{\mbox{Int}}
\newcommand{\Bo}{\Box^{n}_{i}}
\newcommand{\Di}{{\cal D}}
\newcommand{\gd}{{\underline \gamma}}
\newcommand{\gu}{{\underline g }}
\newcommand{\ce}{\mbox{III}}
\newcommand{\be}{\mbox{II}}
\newcommand{\M}{\cal{M}}
\newcommand{\F}{\cal{F_{\eta}}}
\newcommand{\Ci}{\bf{C}}
\newcommand{\ai}{\mbox{I}}
\newcommand{\dupap}{\partial^{+}}
\newcommand{\dm}{\partial^{-}}
\newenvironment{note}{\begin{sc}{\bf Note}}{\end{sc}}
\newenvironment{notes}{\begin{sc}{\bf Notes}\ \par\begin{enumerate}}%
{\end{enumerate}\end{sc}}
\newenvironment{sol}
{{\bf Solution:}\newline}{\begin{flushright}
{\bf QED}\end{flushright}}

\title{A Markov partition for the Feigenbaum dynamics}
\date{}

\author{Genadi Levin
\\
\small{Inst.\ of Math., Hebrew University, Jerusalem 91904, Israel}\\
\small{levin@math.huji.ac.il}\\
}
\normalsize
\maketitle
\abstract{We build a Markov partition of the
Feigenbaum dynamical plane, and a system
of ``external rays'', which are invariant
under rescaling. Applications to geometry
of the Feigenbaum map are given.}

\

\section{Introduction and main results}
``Puzzle" is a Markov partition for the dynamics of a 
polynomial map $z\mapsto z^r+c$. It was introduced by J.-C. Yoccoz 
for non-renormalizable quadratic polynomials
to study the problem of local connectivity 
of Julia sets and the Mandelbrot set~\cite{H}.

We present here a different construction
to build a countable Markov partition
for the famous Feigenbaum map~\cite{F},~\cite{EL},~\cite{E}. 
Remind that this is 
an infinitely renormalizable map. 
The key point is that our partition is invariant
under the rescaling.
To obtain this property, we first construct
an ``external ray'', or access (i.e. a simple arc
disjoint with the Julia set of the map) to the critical point, which
is invariant under the rescaling. Moreover, such ``ray''
appears to be unique, it is a fractal and a quasicircle, see Theorem A.
Taking preimages of this ``ray''
we obtain the Markov partition; 
it is described in Theorem B (see also
Theorem C (a)). Roughly speaking, it says that there exists
(and, in fact, unique) 
a ``shape'', so that its images
under the action of univalent maps with uniformly
bounded distortions (which are inverse branches and rescalings)
tile the plane, and so that the boundary
of any ``tile'' touches a (rescaled) Julia set
at a unique point, which is a (rescaled)
preimage of the critical point. Then we show (Theorem C)
that there is a system of ``external rays'', which is
invariant under the rescaling. 
For other applications, see Theorem C and Corollary.

\

The Feigenbaum map $g$ is a unimodal interval map
with a given criticality $r>1$
at the critical point and such that $g$ is
a fixed point of the period-doubling renormalization
operator~\cite{F}. In other words, it is a solution of 
the Cvitanovic-Feigenbaum equation
\begin{equation}\label{feigeq}
g(x)=\alpha g(g(x/\alpha)),
\end{equation}
where $\alpha<-1$ is the Feigenbaum universal constant
(rescaling in the dynamical plane).
The Feigenbaum map exists for every $r>1$~\cite{La},~\cite{E1}.
The uniqueness is proved for every $r$ even integer~\cite{S}.
Moreover, in this case the map is universal~\cite{F},~\cite{S},~\cite{McM}.
It means that the sequence of renormalizations
of any smooth unimodal map
(with the Feigenbaum combinatorics and fixed criticality $r$)
converges to a map which depends on $r$ only.
This is the Feigenbaum map.

Let $g$ be the Feigenbaum map with the
criticality $r$ of the critical point,
where $r>1$ is a fixed even integer. 
As usual, $g$ is normalized so that
it is an even unimodal map of $[-1, 1]$ into itself,
with the critical point at $0$ and the critical value $g(0)=1$.
It is proved~\cite{EL},~\cite{E}  
that $g$ extends to an analytic polynomial-like map~\cite{DH}
$$g:\Omega_0\to \Omega$$
with the single critical point at $0$,
where $\Omega={\C}\setminus {\R}\cup (\alpha, \alpha^2)$ 
is a slit complex plane, $\overline {\Omega_0}\subset \Omega$,
and $\Omega_0$ {\it is a bounded topological disc};
moreover, if one denotes
$$\Pi=\{z: \arg(z)\in (0, \pi/r)\},$$
then the boundary of $\Omega_0$
contains exactly one {\it singular} point $c$ in $\Pi$ 
such that the limit values of $g$ at $c$
is infinity.

The dynamical system $g:\Omega_0\to \Omega$ is the subject
of the paper.
(So that whenever we speak about iterates $g^n$ or $g^{-n}$,
$n\ge 0$, we mean iterates of the above map.)

There is a representation $g(z)=F(z^r)$
where $F$ is analytic in an appropriate domain.
In particular,
$\Omega_0$ is symmetric w.r.t. the rotation by
the angle $2\pi/r$,
as well as w.r.t. the complex conjugation
$$z\mapsto z^\ast.$$ 

Denote by $J$ the Julia set of $g$:
$J=\{z: g^n(z)\in \Omega_0, n=0,1,2,...\}$.
As it is proved in~\cite{McM},

{\it the union 
$$J_\infty=\cup_{n=0}^\infty \alpha^n J$$
of increasing sequence of compacts
$\alpha^n J$, $n\ge 0$, is dense in the plane.}

This is the reason of uniqueness of the curves
constructed in Theorem A and Theorem C (b)-(b').
Note that the dense set
$J_\infty$ is connected (but not closed).

By a curve, or arc we usually mean an {\it open} one,
i.e. without the end points.
A curve with the end points is called
{\it closed} (this will be clear from the context).

Denote by $HD(E)$ the Hausdorff dimension of
a subset $E$ of the plane.

\

{\bf Theorem A (invariant access at zero)}
{\it There exists a unique simple unbounded curve $L$
which lies in $\Pi$
and terminates at the critical point $0$, and
which has the following properties {\bf (1)-(2)}:

{\bf (1)} $L$ is invariant under the rescaling
$z\mapsto |\alpha|\cdot z$,

{\bf (2)} $L$ is disjoint with the Julia set
$J$ of the polynomial-like map $g$; in fact, 
$$L\cap J_\infty=\emptyset.$$

Additionally, $L$ satisfies the following properties:

{\bf (3)} $1<HD(L)<2$
and $L$ has no tangent at any point.
Moreover, $L$ is a quasicircle (i.e., a quasiconformal image
of the real line),

{\bf (4)} $c\in L$;
if one denotes by $L_n$ the arc of $L$ between
the points $c/|\alpha|^n$ and  $c/|\alpha|^{n+1}$,
$n=0,1,2,...$, and by $L_\infty$ the infinite arc
of $L$ between $c$ and $\infty$, then, for $n=0,1,...,$: 
\begin{itemize}
\item 
(4a) $L_\infty\subset \Omega\setminus \Omega_0$, 
\item 
(4b) the arc
$L_n$ lies in the ``annulus'' 
$g^{-2^n}(\Omega)\setminus g^{-2^{n+1}}(\Omega)$
and joins points $c/|\alpha|^n\in \partial g^{-2^n}(\Omega)$ 
and  $c/|\alpha|^{n+1}\in \partial g^{-2^{n+1}}(\Omega)$, 
\item 
(4c) the curve $g^{2^n}({L_n})$  
is equal to either $(L^{(n)})^\ast$ or
$-(L^{(n)})^\ast$, where $L^{(n)}=L_\infty/|\alpha|^n$ 
is the (infinite) arc of $L$ between
$c/|\alpha|^n$ and $\infty$,
\end{itemize}

{\bf (5)} the intersection of the curve
$g^{2^n}(L\cap g^{-2^n}(\Omega))$ with either $L^\ast$
or $-L^\ast$ is equal to either $(L^{(n)})^\ast$ or
$-(L^{(n)})^\ast$; moreover,
for any $k=2^{n_1}+2^{n_2}+...+2^{n_l}$ 
with $n_1>n_2>...>n_l$, the intersection of
the curve $g^k(L\cap g^{-k}(\Omega))$
with one of the curves $\pm L, \pm L^\ast$ is equal to one 
of the infinite arcs among 
$\pm L_\infty/|\alpha|^{n_l},
\pm L_\infty^\ast/|\alpha|^{n_l}$.}

\

Considering all the preimages of $L\cup L^\ast$, we build a Markov 
partition (even tiling, see below) of $g$,
which is invariant under the rescaling. Let us define it precisely.
Consider two infinite open ``sectors'' $R_0, R_0^\ast$, where
$R_0$ is bounded by the curves $L$ and $-L^\ast$. 
Let us rotate the curves $L$ and $L^\ast$ 
by the angles $2\pi j/r, 0\le j\le r-1$. Then we get $2r$
curves, which
divide either 
``sector'' $R_0, R_0^\ast$
into $r-1$ smaller ``sectors'' $R_{0,1}, R_{0,2},..., R_{0,r-1}$
(in counterclockwise direction)
and $R_{0,1}^\ast, R_{0,2}^\ast,..., R_{0,r-1}^\ast$ respectively.
They are said to be the set ${\M}_0$ of $2(r-1)$ pieces of the 
{\bf depth} $0$.
Now, for every $n=1,2,...$, 
we define the set ${\M}_n$ of pieces of {\bf depth} 
$n$ as the set of components
of $g^{-n}(R)$, $R\in {\M}_0$. Observe that the set of $1$-depth pieces
consists of $2r(r-1)$ bounded Jordan domains.
Hence, all pieces of all depths($\ge 1$) are
bounded Jordan domains. Closed pieces are said to be
the closures of the pieces. For a piece $R$ of depth $n\ge 1$, define
its base points'' $x_R, x_R^\infty\in \partial R$ as follows:
if $g^n(R)\in {\M}_0$ then $g^n(x_R)=0$, $g^n(x_R^\infty)=\infty$.
For example, for the pieces of depth $1$, their base points 
are rotated (by $2\pi j/r, j\in {\Z}$)
and conjugated pairs of points to the
pair $x_0, c$,
where $x_0>0$ is such that $g(x_0)=0$.

\

{\bf Theorem B (the Markov partition)}
{\it The set ${\M}=\cup_{n\ge 1} {\M}_n$ of all pieces satisfies the
following properties.

{\bf (a) (invariance)} if $R\in {\M}_n$, then $R/\alpha\in {\M}_{2n}$.
Besides, $\partial R\cap J_\infty=\{x_R\}$.

{\bf (b) (Markov property)} any two pieces $R, R'$ are either disjoint or 
one contains the other one; moreover, $\partial R\cap \partial R'$
is either (i) empty, or (ii) a non-trivial arc,
which contains one of the base points
$x_R^\infty$, $x_{R'}^\infty$, or (iii) $R\cap R'=\emptyset$
and $\partial R\cap \partial R'=\{x_R\}=\{x_{R'}\}$.

{\bf (c) (tiling)} the union of all closed pieces
together with the real line rotated by the angles
$2\pi j/r, j=1,...,r-1,$ covers a neighborhood
of zero.

{\bf (d) (bounded distortion)} all pieces are roughly
of the same shape; more precisely, for any piece $R$ there exists
a minimal $n\ge 0$ such that: 
\begin{itemize}
\item
(d1) $g^n(R)$ is one of the pieces
of the form $G/\alpha^k$, for some piece $G$ 
of depth one and for some $k\ge 0$,
\item
(d2) the map $g^n: R\to G/\alpha^k$ extends
to a univalent map onto $V(G)/\alpha^k$,
where $V(G)$ is a neighborhood of $\overline G$,
so that $V(G)$ is fixed for each of $2r(r-1)$
closed pieces of depth one.
\end{itemize}
Moreover, there exist $C>0, 0<\lambda<1$, such that,
for every $m$ different pieces $R_m\subset R_{m-1}\subset...\subset R_1$
one has diam$(R_m)\le C\lambda^m$diam$(R_1)$.
}

\
 
Introduce the extended set
of pieces ${\M}^{ext}=\cup_{n\in \Z} \{\alpha^n R: R\in {\M}\}$.
The elements of ${\M}^{ext}$ are also called pieces. Theorem B (a), (c)
shows that ${\M}\subset {\M}^{ext}$
and {\it the set of closures of all (pairwise disjoint) pieces of ${\M}^{ext}$
(united with the rotated real axis) forms a tiling of the plane.}

Let us arrange the set of all pieces ${\M}^{ext}$ by inclusion.
Then there is a countable number of the pieces
of the {\bf first level} (so that they are not containing
in bigger pieces), the {\bf second level} etc:
every piece of {\bf level $n$} is contained in some pieces
of levels up to $n$ and contains pieces of bigger levels. 

Let us introduce the following sets:

$B=\cup_{n\ge 0} \alpha^n B_0$,
where
$B_0$ is the set of all $z$ such that, for some $n\ge 0$,
$g^n(z)=0$.


$Y=\cup_{n\ge 0} \alpha^n Y_0$,
where
$Y_0$ is the set of all $z$ such that, for some $n\ge 0$,
$g^n(z)$ lies in the real axis. In particular,
$B_0\subset Y_0\subset J$.

$\partial {\M}^{ext}$ is said to be the set of
all $z$ such that $z\in \partial R$ for some $R\in {\M}^{ext}$.
In particular, $B\subset \partial {\M}^{ext}$.

Finally, $J^{ess}={\C}\setminus (Y\cup \partial {\M}^{ext})$.
Note that $J\setminus Y_0\subset J^{ess}$.

\

We use Theorems A-B to prove

\

{\bf Theorem C (geometry of the Feigenbaum plane)}

{\it

{\bf (a)} for every 
$z\in J^{ess}\cup\partial {\M}^{ext}\setminus Y=
{\C}\setminus Y$ there is a sequence 
of closed nested pieces $(\overline{R_n(z)})_{n\ge 0}$ where 
$R_n(z)\in {\M}^{ext}$
is a piece of level $n$ such that $z\in \overline{R_n(z)}$;
moreover (see Theorem B (d)),
diam$R_{n}(z)\le C\lambda^n$diam$R_1(z)$, with the universal
$C>0, \lambda<1$. If $z\in J\setminus Y_0$, then all $R_n(z)\in {\M}$;

{\bf (b)}  let $x, y\in \overline \Pi$.
There exists a simple curve $l(x, y)$, such that:

(b1) it joins the point $x$ and $y$,

(b2) it lies in $\Pi\setminus J_\infty$.

Moreover, $l(x, y)$ consists of arcs of the boundaries of pieces
of ${\M}^{ext}$,
and if $x, y$ belong
to a closed piece, then $l(x, y)$ 
lies in the same closed piece;
finally, if 
$x, y\in \Pi\setminus Y$, then $l(x, y)$ is the
unique simple curve that satisfies (b1)-(b2);

{\bf (b')} given $x\in \overline \Pi$,
define an ``external ray''
$l\subset \Pi$ to the point $x$ as the curve $l(x, \infty)$
(i.e. the simple curve that joins $x$ and $\infty$
and is contained inside the curve $l(x, c)$
extended by the infinite arc $L_\infty$).
Then: (i) for any $n\in {\Z}$, $l(|\alpha|^n x, \infty)=
|\alpha|^n l(x, \infty)$, and (ii)
for any two ``external rays''
$l_1, l_2$, the intersection
$l_1\cap l_2$ is equal to a simple curve joining
a point $y\in l_1\cap l_2$ and infinity, where
$y$ is such that, for some non-negative $n, k$,
the point $g^n(y/\alpha^k)$
belongs to the set $\{ c\exp(i2\pi j/r), c^\ast\exp(i2\pi j/r),
j=0,...,r-1\}$,

{\bf (c)} there exists $C$ such that,
for every point $z\in \overline \Pi$, 
there exists a rectifiable curve inside of $J_\infty\cap \overline{\Pi}$
with the length less than $C\cdot |z|$, 
which connects $z$ and $0$.
}

\

We complete the system of ``external rays'' 
from $\Pi$ to the whole plane in a natural way:

\

{\bf Complement to (b').} {\it Applying to the
``external rays'' in $\Pi$ complex conjugation and rotations
by the angles $2\pi j/r$ ($1\le j\le r-1$),
we define the set of ``external rays'' in the whole plane
(in particular, to any point of the plane).
Specifically, there are precisely $2r$ ``rays''
to the critical point $0$ (obtained as above
from the ``ray'' $L\subset \Pi$). Hence, one can define
precisely $2r$ ``rays'' to any point
of the set $B$. Moreover, if $x\in {\R}\setminus B$, 
then, as it easy to see, there are
precisely $2$ ``rays'' to $x$
(at least two, by the complex conjugation, and at most two,
by topological considerations). Therefore, one can define
precisely $2$ ``rays'' to any point of $Y\setminus B$.
To any other point of the plane,
the number of ``rays'' is $1$ (by Theorem C (b)).}

\

It is known~\cite{E},
that $g$ has a maximal analytic continuation to
a domain $W$,
which is simply connected. Moreover, 
$W=\cup_{n\ge 0} \alpha^n g^{-2^n}(\Omega)$
~\cite{E},~\cite{B}, and $W$ is dense
in the plane~\cite{McM} (because the dense
set $J_\infty\subset W$).
Note that $W$ is invariant under the complex conjugation
and rotation by $2\pi/r$.
Then Theorems A (4b) and C (b) give us 

\

\begin{coro}\label{anal}

{\bf (a)} $L_\infty\subset \partial W$,

{\bf (b)} $\partial W$ is a dendrite, that is: 
\begin{itemize}
\item 
(b1) for every two points
$x, y\in \partial W\cap \Pi$ there exists a unique
curve, namely, $l(x, y)$, see Theorem C (b),
that joins $x$ and $y$ inside $\partial W\cap \Pi$,
\item
(b2) if $x, y\to z\in \partial W$, then 
the diameter of $l(x, y)$ tends to zero.
\end{itemize}
\end{coro}
\begin{proof} Part (a) follows immediately from
Theorem A (4b). To prove (b1), assume the contrary, i.e.,
assume that, for some $x\not=y$ in $\partial W\cap \Pi$,
the curve $l(x, y)$ is not contained in $\partial W$.
Since $\partial W$ is connected (on the Riemann sphere),
then $l(x, y)\cup \partial W$ separates the sphere.
On the other hand, both sets $\partial W$ and $\overline{l(x, y)}$
are disjoint with the dense and
connected set $J_\infty$, a contradiction. 
By the same reason, the curve $l(x, y)$ is the unique
simple curve in the plane connecting $x, y$ inside of $\partial W$.

For the proof of (b2), see Sect.~\ref{theoc}.
\end{proof}

\section{Proof of Theorem A}
We construct the curve $L$ as the union 
$\cup_{n\in \Z} |\alpha|^n I$,
where $I$ is a simple (closed, i.e. with the end points) curve 
joining $c$ and $c/|\alpha|$. In turn, the curve $I$
is going to be the limit set
of a 
{\bf finite} system of conformal maps (iterated function system).

\subsection{Some basic facts}
Denote by
$\h^\pm$ the upper and lower half-planes. 
Let also $x_0$ be the unique point $0<x_0<1$,
such that $g(x_0)=0$.
The following facts {\bf 1-3} (see~\cite{EL},~\cite{E}) will be of crucial
importance for the proof.

{\bf 1.} Let
$u$ be the branch of $g^{-1}$ such that $u(1)=0$ and
$u(1/\alpha)=1$. Then $u$ extends as a univalent 
(symmetric under the complex conjugacy) function
to the slit plane $\h^+\cup \h^-\cup (\alpha, 1)$, and
$u(\h^-)\subset \Pi$.
It satisfies the equation 
$$u(z)=-\alpha u(u(z/\alpha)), z\in \h^+\cup \h^-\cup (\alpha, 1).$$

{\bf 2.} Introduce the following antiholomorphic maps of $\h^+$ into
itself:
$$u^\ast: z\mapsto u(z^\ast),$$
$$\chi:z\mapsto |\alpha| u^\ast(z),$$ 
Then 
$u^\ast$, $\chi$ extend to homeomorphisms 
of $\overline \h^+$ into itself, and
$\chi: \h^+\to \h^+$
has a unique (globally attracting) fixed point, namely, the
singular point $c\in \Pi$.

{\bf 3.} Denote
\begin{equation}\label{omega}
\omega=u^\ast(\h^+).
\end{equation}
Then $\omega$ is a Jordan domain.
Its boundary consists of 
a real segment
$[0, |\alpha|x_0]=u^\ast([\alpha, 1])$, a sequence
of (closed) analytic arcs $\tau_n, n\ge 0$, and the point $c$.
Here $\tau_0$ is the segment $\exp(i\pi/r)[0, |u(\alpha^2)|]=
u^\ast([1, \alpha^2])$,
and $\tau_{n+1}=\chi(\tau_n)=u^\ast((-1)^{n+1}[|\alpha|^{n+1},
|\alpha|^{n+3}])$, $n\ge 0$.
In particular, $\{\tau_n\}$ converge to the point $c$
exponentially fast.
For every $n$, $\tau_n, \tau_{n+2}$ 
meet at a joint end point at the angle $\pi/r$.
 
\

Observe (after~\cite{EL}) that the map $u^\ast$
conjugates $z\mapsto \alpha z^\ast$ and $\chi$:
\begin{equation}\label{conj}
\chi\circ u^\ast (z)=u^\ast(\alpha z^\ast).
\end{equation}
Hence, 
\begin{equation}\label{conj2}
\chi^2\circ u^\ast (z)=u^\ast(\alpha^2 z).
\end{equation}
In particular, we have
\begin{com}\label{lens}
If $\gamma$ is any straight ray 
in $\h^+$ starting at zero,
then $u^\ast(\gamma)$ is an (open) analytic arc in $\omega$
joining the points $c$ and $x_0$ in $\partial \omega$,
which has a finite length (because it is invariant under
$\chi^2$). Moreover, since the map $u^\ast$
extends to a univalent map in a neighborhood
of zero,  $u^\ast(\gamma)$ extends to an analytic curve
through the point $x_0$.
\end{com}
 
\subsection{Two iterated function systems}
Note that
$$u^\ast(\omega)=u(u(\h^+))=u(u(\h^-/\alpha))=
u(\h^-)/(-\alpha)=\omega/|\alpha|.$$
Since $\omega\subset \Pi$, we have 
$\omega/|\alpha|=u^\ast(\omega)\subset u^\ast(\h^+)=\omega$
and $u^\ast(\omega)\subset u^\ast(\Pi)$.





Introduce a compact set 
$$X=\overline{u^\ast(\Pi)\setminus u^\ast(\omega)}=
\overline{u^\ast(\Pi)\setminus (\omega/|\alpha|)},$$
and the following {\it two} 
families $S_\infty=\{\psi_k\}_{k=1}^\infty$,
$S=\{\varphi_1,\varphi_2, \varphi_3\}$ 
of conformal maps 
from $\h^+$ into itself defined as follows:
\begin{equation}\label{psidef}
\psi_k(z)=u^\ast(|\alpha|^k z), k=1,2,3,4,...,
\end{equation}
and
\begin{equation}\label{phidef}
\varphi_1=\psi_1, \varphi_2=\psi_2, \varphi_3=\chi^2. 
\end{equation}
The eq.~(\ref{conj}) implies
\begin{equation}\label{psi}
\psi_{2j}(z)=\chi^{2j}(u^\ast(z)), (j\ge 1) \ \
\psi_{2j+1}(z)=\chi^{2j+1}(u^\ast(-z^\ast)), (j\ge 0).
\end{equation}
Therefore,
\begin{equation}\label{ind}
\psi_{2j+2}=\varphi_3^{j}\circ \varphi_2, \ \
\psi_{2j+1}(z)=\varphi_3^{j}\circ \varphi_1, \ \ (j\ge 0).
\end{equation}
\begin{com}\label{gen}
The latter means that the infinite system $S_\infty$ 
is {\it induced} by the finite system $S$.
So, it is enough to study the system $S$ (see also Comment~\ref{s,sinfty}).
Although the system $S_\infty$ arises naturally if one wants to 
construct the curve $L$, the system $S$ is simpler to deal with.
In the sequel,
$S_\infty$ is more convenient
for the study of
properties related to rescaling
((1)-(2), (4)-(5) of Theorem A)
while $S$ is used to prove the metric
properties (3), Theorem A.
\end{com}
Denote by $\rho$ the hyperbolic metric of $\h^+$.
Note that maps $u^\ast$ and $\chi$ as conformal maps
of $\h^+$ into itself are contractions
w.r.t. the metric $\rho$.

\begin{lem}\label{ifs}
The family $S$ obeys the following properties.

{\bf (1)} $\Int(X)$ is a Jordan domain. The boundary
of $X$ consists of: (i) a sequence of (closed) curves
$\tau_n/|\alpha|$, $n\ge 2$, which converges
to the points $c/|\alpha|\in \partial X$,
so that $\cup_{n\ge 2} \tau_{n}/|\alpha|\cup \{c/|\alpha|\}$ 
is a (closed) curve, which we call the ``bottom'' $B_X$ of $X$;
(ii) the curve $\tau_0\setminus (\tau_0/|\alpha|)$
and a sequence 
of curves $\tau_{2j}$, $j\ge 1$, which converges to the
point $c\in \partial X$,
so that 
$\tau_0\setminus (\tau_0/|\alpha|)\cup \cup_{j\ge 1} \tau_{2j}\cup \{c\}$ 
is a curve, which we call the ``left side'' $L_X$ of $X$;
(iii) the curve $u^\ast(\exp(i\pi/r) \R^+\setminus \tau_0)$,
which we call the ``right side'' $R_X$ of $X$.
Finally, the compact $X$ is contained in $\h^+$.

{\bf (2)} $|\alpha| X\cap X=\cup_{j\ge 1}\tau_{2j} \cup \{c\}$.

{\bf (3)} 
Each $\ph_i$ 
is either holomorphic
or antiholomorphic injection of $\h^+$
into itself, which maps the compact $X$ into itself;
moreover, there exist a domain $V$ and a number $\lambda<1$, 
such that
$X\subset V$, $\overline V\subset \h^+$,
$\ph_i(V)\subset V$, and
$$\rho(\ph_i(x), \ph_i(y))\le \lambda \rho(x, y),$$
for $i=1,2,3$ and all $x, y\in V$.

{\bf (4)} $\varphi_i(\Int(X))\cap \varphi_j(\Int(X))=\emptyset$
for $1\le i<j\le 3$.

{\bf (5)} $\varphi_i(X)\cap \varphi_{j}(X)\not=\emptyset$ if and only if
$|i-j|\le 1$.

{\bf (6)} for any finite word $\bar \e=\{\e_1,...\e_n\}$ of 
$n$ symbols
$\e_i\in \{1,2,3\}$, denote
$\ph_{\bar \e}=\ph_{\e_1}\circ \ph_{\e_2}\circ...\circ \ph_{\e_n}$.
Then $\ph_{\bar \e}$ maps $V$ into itself with uniformly bounded
distortion:
there is $K$ such that
$|\ph_{\bar \e}'(x)|\le K |\ph_{\bar \e}'(y)|$ for every finite
word $\bar \e$ and any $x, y\in V$.
\end{lem}
\begin{proof}
Parts (1)-(2) follow from the definition of $X$ and from the
description of $\omega$.
 
Part (3). For every $k>0$, $\psi_k(X)=u^\ast(|\alpha|^k X)\subset
u^\ast(\overline{\Pi\setminus \omega})=X$.
In particular, $\varphi_i(X)\subset X$, $i=1,2$.
As for $\varphi_3=\chi^2$, one can
use that $\psi_2=\chi^2\circ u^\ast$ and write
$\chi^2(X)=\chi^2\circ u^\ast(\overline{\Pi\setminus \omega})=
\psi_2(\overline{\Pi\setminus \omega})=
u^\ast(\alpha^2\overline{\Pi\setminus \omega})\subset
u^\ast(\overline{\Pi\setminus \omega})=X$.

Let's construct the domain $V$
(the existence of $\lambda$ then follows because
$\overline V\subset \h^+$).
We start with the following general

{\bf Claim.}
{\it For a point $a\in {\bf C}$, a simply connected
domain $U\ni a$, and a map $f$ which
is a holomorphic injection 
of $U^\ast=\{z^\ast| z\in U\}$ into $U$, 
define an antiholomorphic map $f^\ast(z)=f(z^\ast)$, $f^\ast:U\to U$; 
assume that
$f^\ast(a)=a$
and $0<|\rho|<1$, where $\rho=f'(a^\ast)$.
Then there is a holomorphic change
of coordinates $K:U\to \C$, $K(a)=0$, such that
$K\circ f^\ast\circ K^{-1}: w\mapsto \rho w^\ast$.
Moreover, if $Y\subset U$ is compact, 
such that $f^\ast(Y)\subset Y$, then
there is an arbitrary small 
neighborhood $W_Y$ of $Y$, such that
$f^\ast(W_Y)$ is compactly contained in $W_Y$.}

It's enough to prove the existence of $K$ locally,
around $a$. Notice that the second iterate $(f^\ast)^2$
of $f^\ast$ is holomorphic around its fixed point $a$,
and its multiplier $((f^\ast)^2)'(a)$ is equal to $|\rho|^2$.
Hence, there is a holomorphic change $K$, $K(a)=0$, $K'(a)=1$,
such that $K\circ (f^\ast)^2\circ K^{-1}: w\mapsto |\rho|^2 w$.
Then the map $K\circ f^\ast\circ K^{-1}$
is antiholomorphic around $a$, hence, there is
a univalent at $w=0$ map $F$, such that
$K\circ f^\ast\circ K^{-1}=F^\ast$.
If $F(w)=\rho w + a_n w^n+...$, then
one sees by substitution that $a_n=0$.
Thus $F(w)=\rho w$, which proves that 
$K\circ f^\ast\circ K^{-1}(w)=\rho w^\ast$.
Now, for every $\epsilon>0$, the neighborhood 
$W_Y$ can be defined as
$K^{-1}(W^\epsilon)$, where
$W^\epsilon$ is the $\epsilon$-neighborhood
of $K(Y)$.  

Let's apply Claim to the map
$\varphi_1$ and the point $c/|\alpha|$. It is relevant, because
$\varphi_1=\psi_1=|\alpha|^{-1}\circ \chi\circ |\alpha|$,
hence, $c/|\alpha|$ is the (global) attracting fixed point
of the map $\varphi_1:\h^+\to \h^+$.
Besides, $c/|\alpha|$ lies in the ``bottom''
$B_X$ of $X$. Note that $\varphi_1(B_X)\subset B_X$.
Hence, there is a small enough
neighborhood $W_B$ of $B_X$, such that
$\varphi_1(W_B)$ is compactly contained in $W_B$.

Now note that the ``left side'' $L_X$ of $X$ as well as 
its ``right side'' $R_X$ is invariant under the holomorphic map
$\ph_3$. Hence, one can choose a small enough neighborhood
$W_L$ of $L_X$, such that $\ph_3(W_L)$ is compactly
contained in $W_L$. Moreover, note the following (a)-(c):

(a) $\ph_1$ maps the ``left side'' $L_X$ onto a curve, which is
an arc in ``right side'' $R_X$ extended by an arc in $Int(X)$,
and it maps $R_X$ into $Int(X)$ except for an end point in $B_X$;
(b) $\ph_2$ maps $L_X$ onto a curve, which is
an arc in $R_X$ extended by an arc in $Int(X)$, and it maps
$R_X\cup B_X$
into $Int(X)$  except for an end point in $B_X$; 
(c) $\ph_3$ maps $B_X$ into $Int(X)$ except for
end points in $L_X\cup R_X$. 

It follows that there is a small neighborhood
$W_R$ of $R_X$, such that $\ph_3(W_R)$ is compactly contained
in $W_R$ and  
$\ph_1(W_R)\cup \ph_2(W_R)\subset X\cup W_B$, moreover, 
$W_L$ is chosen so small 
that $\ph_1(W_L)\cup \ph_2(W_L)\subset W_R$.
Thus one can set $V=X\cup W_R\cup W_L\cup W_B$.
We have proved Part (3). This implies also, by Koebe distortion theorem,
the bounded distortion property (6).

Let's prove parts (4)-(5). For $i=1, j=2$, this 
follows directly from the definition of $\ph_k$.
For $i=2, j=3$, we have
$\ph_2(X)=\ph_3\circ u^\ast(X)\subset \ph_3\circ u^\ast(\omega)
=\ph_3(\omega/|\alpha|)$. Since $Int(\omega/|\alpha|)$
is disjoint with $Int(X)$, then
$\ph_2(Int(X))\subset \ph_3(Int(\omega/|\alpha|))$
is disjoint with $\ph_3(Int(X))$.
On the other hand, $\ph_2(X)\cap \ph_3(X)$ contains
$\ph_2\cap \ph_3\circ \ph_1(X)=\psi_2(X)\cap \psi_3(X)$,
which is a non-empty set, by the definition
of $\psi_k$ and by (2). 
Finally, for $i=1, j=3$, $\ph_1(X)=u^\ast(|\alpha| X)
\subset u^\ast(|\alpha|\omega)$ while
$\ph_3(X)=\ph_3\circ u^\ast(\overline{\Pi\setminus \omega})$,
hence, by~(\ref{conj2}), $\ph_3(X)=
u^\ast(\alpha^2 \overline{\Pi\setminus \omega})$, and 
this set
is disjoint with  $u^\ast(|\alpha|\omega)$
because $|\alpha|\omega$ and $\alpha^2 \overline{\Pi\setminus \omega}$
are disjoint.
\end{proof} 
\subsection{The limit set}
For any infinite word
of symbols $\bar \e=\{\e_n\}_{n>0}$, $\e_n\in \{1,2,3\}$, denote 
$\bar \e|n=\{\e_1,\e_2,...,\e_n\}$
and define the limit set of $S$ as
$$I=\cup_{\bar \e}\cap_{n=1}^\infty \ph_{\bar \e|n}(X).$$
Since $S$ is finite, $I$ is a closed set. 
It is easy to see that
$$I=\cap_{n>0} X_n.$$
Here
$X_n$ is a decreasing sequence of compacts
$X_n=
\cup_{|\bar \e|=n}\ph_{\bar \e}(X)$
where $|\bar \e|=n$ means that the number of symbols
in the finite word $\bar \e$ is $n$. 
Since $\ph_3(c)=c$ and $\ph_1(c/|\alpha|)=c/|\alpha|$,
then $c, c/|\alpha|\in I$.




\begin{com}\label{s,sinfty}
Let's consider the closure $\overline I_\infty$
of the limit set $I_\infty$ of the infinite system
$S_\infty=\{\psi_k:X\to X\}_{k\in \N}$, where
$I_\infty=\cup_{\bar w} \cap_{n=1}^\infty \psi_{\bar w|n}(X)$,
and $\bar w$ are the infinite words of symbols
$w_n\in \N$. Then it is easy to see that
$I=\overline I_\infty$. Indeed,
since $S_\infty$ is induced by $S$,
$\overline I_\infty\subset I$. To prove the opposite
inclusion, notice that $\ph_3^n(X)\to c$ ($n\to \infty$),
and
consider any point
$a\in I$, $a=\cap_{n=1}^\infty \ph_{\e|n}(X)$,
where $\e$ is an infinite word of symbols
$\e_n\in \{1,2,3\}$. If $\e_n=3$ for all $n>n_0$, 
and $\e_{n_0}\in \{1,2\}$, then
$a=\ph_{\e|n_0}(c)$ and $\ph_{\e|n_0}=\psi_{w}$
for a finite word $w$ (in the alphabet $\N$ of $S_\infty$).
If $\e_{n_j}\in \{1,2\}$, for $n_j\to \infty$, then 
there is an infinite word $\bar w=\{w_n\}_{n=1}^\infty$,
$w_n\in \N$, such that
$\ph_{\e|n_j}=\psi_{\bar w|j}$, and hence
$a\in I_\infty$.
\end{com}

Our first aim is to show that $I$ is a simple (closed, i.e. with
the end points $c, c/|\alpha|$)
curve. Since $X\setminus \cup_{1\le i\le 3}\ph_i(X)$ 
contains an open non-empty disc, and because
of contraction, the Lebesgue measure (area)
of $I$ is zero. In particular,
$I$ has no interior. Clearly, $I$
is also locally connected. 
To be more precise, if
$x\in \cup_{\bar \e}\cap_{n=1}^\infty \ph_{\bar \e|n}(Int(X))$,
there is a unique word $\bar \tau$, such that
$\ph_{\bar \tau|n}(I), n\ge 1$
form a sequence of connected closed neighborhoods
of $x$ shrinking to $x$.
Let now $x\in  \ph_{\tau}(\partial X)\cap I$,
for a finite word $\tau$, $|\tau|=n$.
Assume first that $n=1$.
As it follows from (a), (b), (c) (see the proof of Lemma~\ref{ifs}), 
$X_n\cap \partial X\subset \partial \ph_1^{n-1}(X)\cup \partial \ph_3^{n-1}(X)$
for $n\ge 2$,
and that $\ph_1^n(X)\to c/|\alpha|$, 
$\ph_3^n(X)\to c$ as $n\to \infty$.
Therefore, 
\begin{equation}\label{twopoints}
I\cap \partial X=\{c/|\alpha|, c\} 
\end{equation}
and, moreover,
$I\cap \cup_{1\le i\le 3} \partial \ph_i(X)$
consists of the following $4$ points $a_j, 1\le j\le 4$,
with corresponding sequences of closed connected
neighborhoods:
\begin{equation}\label{a4a1}
a_4=c=\cap_{n=1}^\infty \ph_3^n(I), \ \
a_1=c/|\alpha|=\cap_{n=1}^\infty \ph_1^n(I), 
\end{equation}
\begin{equation}\label{a2}
a_2=\ph_1(a_4)=\ph_2(a_1)=
\cap_{n=1}^\infty \ph_1(\ph_3^n(I))\cup \ph_2(\ph_1^n(I)), 
\end{equation}
\begin{equation}\label{a3}
a_3=\ph_2(a_4)=\ph_3(a_1)=
\cap_{n=1}^\infty \ph_2(\ph_3^n(I))\cup \ph_3(\ph_1^n(I)).
\end{equation}
If now $n>1$, then $x$ is an image by $\ph_{\tau|n-1}$
of one of the points $a_j, 1\le j\le 4$.
Thus $I$ is locally connected at every point.

We have just proved that $I$ is a dendrite.
In particular, for every two points $x, y\in I$
there exists a unique simple arc $[x,y]$
connecting $x, y$ inside $I$.
Now we show that $I=[c/|\alpha|, c]$.  
Assume there is some $z\in I\setminus [c/|\alpha|, c]$.
Then we get a (branch) point $d\in [c/|\alpha|, c]$, i.e.
such that $[z, c]\cap [c/|\alpha|, c]=[d,c]$
and $[z, c/|\alpha|]\cap [c/|\alpha|, c]=[d,c/|\alpha|]$. 
If $d\in \cup_{\bar \e}\cap_{n=1}^\infty \ph_{\bar \e|n}(Int(X))$,
then
passing to a big scale by some $(\psi_{\bar \e})^{-1}$ 
we get that $I$ has at least $3$ different points with the boundary
of $X$, a contradiction with~(\ref{twopoints}). 
If $d\in  \ph_{\tau}(\partial X)\cap I$,
for a finite word $\tau$, $|\tau|=n$,
then applying $\ph_{\tau|n-1}^{-1}$ one can assume 
from the beginning that the branch point $d$ 
is one of the points $a_j, 1\le j\le 4$,
also that, for some $1\le i\le 3$, 
$[d,z]\cap \ph_i(X)$ contains a non-trivial
subarc $[d,z']$, and either $[c,d]$
or $[c/|\alpha|, d]$ intersects $\partial \ph_i(X)$.
Then $\ph_i^{-1}(d)$ is either $c$ or $c/|\alpha|$.
If, say, it is $c$, then $[c/|\alpha|, c]\cup 
\ph_i^{-1}([d,z'])$ is a subset of $I$,
and it intersets $\partial \ph_3^n(X)$, $n$ big,
in at least $3$ points, a contradiction
with~(\ref{twopoints}).

\subsection{The curve $L$. Proof of (1)-(2) and (4)-(5) of Theorem A}
Define two curves: $L_0=I\setminus \{c, c/|\alpha|\}$ and 
$$L=\cup_{n\in \Z} |\alpha|^n I.$$
Then $L$ is a curve which satisfies the property (1)
of Theorem A. 
$L$ is a simple curve because $c\in I\cap |\alpha|I\subset 
(Int(X)\cup \{c, c/|\alpha\})\cap |\alpha|(Int(X)\cup \{c, c/|\alpha\})=
\{c\}$, by Lemma~\ref{ifs}(2).
By Comment~\ref{s,sinfty} and from the
definition of the maps $\psi_k$, we see that
$$u^\ast(L_\infty)=\cup_{k>0} u^\ast(|\alpha|^k (L_0\cup \{c\})=
\cup_{k>0} \psi_k(L_0\cup \{c\})=L_0,$$
i.e.
$$g(L_0)=L_\infty^\ast.$$
Therefore, using the functional equation,
properties (4a)-(4c) follow easily. In particular, (4c) implies that
$L_n\cap J=\emptyset$ for all $n\ge 0$ thus proving also (2).
Let's prove property (5) and the uniqueness.

The uniqueness: assume that $l$ is a curve in $\Pi$ satisfying (1)-(2)
of Thm. A.
Let $x\in l\cap \partial \Omega_0$. Since
$\partial (\Pi\cap \Omega_0)\setminus J_\infty=\{c\}$,
we mush have $x=c$. Hence also 
$c/|\alpha|^n\in l$, $n\ge 0$. If now $l$ is different from $L$,
then (1)-(2) give us a contradiction, as in the proof
of Corollary~\ref{anal}.  

(5): for any $k>0$,
if $x$ lies in the intersection of $\gamma=g^k(L\cap g^{-k}(\Omega))$
and, say, $L$, then the whole 
arc of $L$ from $x$ to $\infty$ also belongs to $\gamma$,
because otherwise the dense set $\cup_{n\ge 0} \alpha^n J$
would be not connected, see the proof of Corollary~\ref{anal}.
Prove other statements.
First, it is easy to see that
$g(L\cap g^{-1}(\Omega))\cap L^\ast=L^\ast_\infty$
(because $g(L_0)=L^\ast_\infty$ and 
$g^2(L_1)=L^\ast_\infty/\alpha$).
From this and from the functional equation, it follows that,
for any $n\ge 0$, 
$g^{2^n}(L\cap g^{-2^n}(\Omega))$ intersects either $L^\ast$
or $-L^\ast$ along either $(L^{(n)})^\ast$ or
$-(L^{(n)})^\ast$. Let's apply this consequently
for $n=n_1, n_2,..., n_l$
when $k=2^{n_1}+2^{n_2}+...+2^{n_l}$
with $n_1>n_2>...>n_l$. 
Then we get the statement.
\subsection{Metric properties of $L$: proof of (3) of Theorem A}
\paragraph{The curve $L$ is nowhere differentiable.}

It is enough to show
\begin{lem}\label{notang}
The curve $I$ has no tangent at any point.
\end{lem}
\begin{proof} We follow~\cite{br}.
Assume the contrary, i.e.
$I$ has a tangent
at a point $a\in I$.
Since $S$ is finite, there is an infinite word
$\epsilon=\epsilon_1 \epsilon_2 ...\epsilon_n ...$, such that
$a=\cap_{n=1}^\infty \phi_{\epsilon|n}(X)$.
Then the maps $f_n:=\phi_{\epsilon|n}:V\to V$, $n\ge 1$,
converge to the point $a$ uniformly in $V$.
Let $b=\lim_{j\to \infty} b_{n_j}$, where
$b_n=f_n^{-1}(a)$. 
One can further assume (possibly passing to 
a subsequence) that all maps
$f_{n_j}$ are holomorphic
(if they are all antiholomorphic,
replace $f_{n_j}$ by complex conjugate to $f_{n_j}$). 
Define
$F_j=(f_{n_j}-a)/f_{n_j}'(b_{n_j})$.
Since $F_j(b_{n_j})=0$, $F_j'(b_{n_j})=1$, and
$F_j$ are univalent, one can assume
(possibly passing to a subsequence) that
$F_j$ converge uniformly in $V$ to a univalent map
$F$, $F(b)=0, F'(b)=1$.
Then $f_{n_j}(z)-a=\lambda_j (F(z)+\delta_j(z))$,
where $\lambda_j=f_{n_j}'(b_{n_j})\to 0$ as $j\to \infty$,
$\lambda_j\not=0$, and
$\delta_j\to 0$ uniformly in $V$.
Consider any $z\in I$, $z\not=b$.
Then $F(z)\not=0$.
If one fixes a sequence $j_k\to \infty$, such that
$\arg \lambda_{j_k}\to \theta$, then 
$\arg(f_{n_j}(z)-a)\to \theta+\arg F(z)$.
By our assumption, we conclude that
the function $\arg F(z)$ is a constant in $z\in I$.
Therefore, $I$ is an analytic arc.
In particular, it is analytic in a neighborhood
of the point $u^\ast(|\alpha|c)\in L_0=I\setminus \{c, c/|\alpha|\}$.
But $u^\ast(L_\infty)=L_0$.
Therefore, the curve $L=\cup_{n\in \Z} |\alpha|^n I$ 
is analytic in a neighborhood
of the point $c$, and, hence, the curve  $L$ is analytic
at any point.
Consider the curve $u^\ast(L)$. It is also
analytic and joins the points $x_0\notin L$ and
$c\in L$. Moreover, $L_0\subset L\cap u^\ast(L)$. 
Thus $L_0$ has two different analytic continuations,
a contradiction.
\end{proof} 
This proves that $L$ has no tangent at any its point.
\paragraph{Hausdorff dimension of $L$.}

Here we prove that  $HD(I)>1$.
The method of the proof is very similar to
the one for the Julia set of $z^2+\e$~\cite{su3}.
First, we show the existence of an $h$-conformal
measure for the system $S$ (Lemma~\ref{meas}), 
then we derive from here that
$h=HD(I)$ and that the $h$-Hausdorff measure
of $I$ is between $0$ and $\infty$ (Corollary~\ref{frost}).
Hence, if $h=1$, then the curve $I$ is rectifiable, 
therefore,
has tangent at many points, a contradiction with Lemma~\ref{notang}. 
Here we only state Lemma~\ref{meas} and Corollary~\ref{frost}.
The proof is very close to the proof of corresponding
statements in~\cite{sum} and~\cite{MU}, and postponed until Appendix
(one could also apply directly the theory developed in~\cite{MU}).
\begin{lem}\label{meas}
There exist $0<h\le 2$ and a probability 
measure $\mu$ with its support $supp(\mu)=I$,
which is $h$-conformal, i.e., 
$\mu(\ph_i(A))=\int_A |\ph_i'(x)|^h d\mu(x)$ for
$i=1,2,3$ and any Borel set $A\subset X$. 
Moreover, $\mu$ has no atoms.
\end{lem}
\begin{coro}\label{frost}
$h=HD(I)$. Moreover, the $h$-Hausdorff measure of $I$ is finite
and positive.
\end{coro}
As we noticed before, Lemma~\ref{meas}, Corollary~\ref{frost},
and Lemma~\ref{notang} imply that $HD(L)>1$.

\paragraph{$L$ is a quasi-circle.}

The inequality $HD(I)<2$ follows from Lemma~\ref{meas} and
~\cite{MU}, Thm. 4.5. Hence,
$HD(L)=HD(I)<2$.
But we prove a stronger statement: $L$ (or, equivalently, $I$) is a
quasicircle.
According to~\cite{A}, it is enough to prove that
$I$ satisfies $M$-condition: there is $M$ such that
for every $3$ points $\xi_1, \xi_2, \xi_3$ of the curve $I$ such that
$\xi_3$ lies between $\xi_1$ and $\xi_2$,
$|\xi_3-\xi_1|/|\xi_2-\xi_1|\le M$. 
Let's agree that the point $\xi_1$ lies between $c/|\alpha|$
and $\xi_3$.
The $M$-property follows from the following two
observations.

(a) Because of the bounded distortion property Lemma~\ref{ifs}(6),
one can assume that the points $\xi_1$ and $\xi_2$
belong to different $\ph_k(I), \ph_n(I), 1\le k<n\le 3$:
$\xi_1\in \ph_k(I), \xi_2\in \ph_n(I)$.

(b) If $n-k=2$, then there is $\sigma>0$ such that $|\xi_1-\xi_3|\ge \sigma$,
for all $\xi_1, \xi_3$ as above, and
the existence of some $M=M(\sigma)$ follows.
Assume now that $n-k=1$. For example, $\xi_1\in \ph_2(I)$,
$\xi_2\in \ph_3(I)$. Note that $\ph_2(I)$
and $\ph_3(I)$ has in common a single boundary point
$a=a_3=\ph_2(c)=\ph_3(c/|\alpha|)=
\cap_{n=1}^\infty \ph_2(\ph_3^n(X))\cup \ph_3(\ph_1^n(X))$.
By the triangle inequality,
enough to prove the $M$-property for $\xi_3=a$. So, assume $\xi_3=a$.
Note that $\xi_1\in \ph_2(\ph_3^n(X))$, $\xi_2\in \ph_3(\ph_1^m(X))$,
for some $n,m\ge 0$, and the sets $\ph_3^k(X)$, 
$\ph_1^{2k}(X)$ converge (as $k\to \infty$)
to the fixed points $c$, $c/|\alpha|$ 
of $\ph_3$, $\ph_1$ resp.
with the same speed, because 
$\ph_1^2=|\alpha|^{-1}\circ \ph_3\circ |\alpha|$.
Hence,
$|\xi_1-a|/|\xi_1-\xi_2|$ is comparable to
$|\tilde \xi_1-a|/|\tilde \xi_1-\tilde \xi_2|$,
where $\tilde \xi_1\in \ph_2(\ph_3^{n_0}(X)\setminus \ph_3^{n_0+1}(X))$, 
$\tilde \xi_2\in \ph_3(\ph_1^{m_0}(X)\setminus (\ph_1^{m_0+2}(X))$,
with some $n_0, m_0\ge 0$,
and, additionally, either $n_0$ or $m_0$ is zero.
Besides, any two sets of the form
$\ph_2(\ph_3^{n}(X)\setminus \ph_3^{n+1}(X))$, 
$\ph_3(\ph_1^{m}(X)\setminus \ph_1^{m+2}(X))$
are disjoint. All this implies that the distance
between $\tilde \xi_1, \tilde \xi_2$ is bounded from below
by the (positive) distance between the sets
$\ph_2(\ph_3^{n_0}(X)\setminus \ph_3^{n_0+1}(X))$,
$\ph_3(\ph_1^{m_0}(X)\setminus (\ph_1^{m_0+2}(X))$.
The $M$-property follows.
\section{Proof of Theorem B}
{\bf (a):} note that, for any $n\ge 0$, ${\M}_n$ 
is invariant w.r.t. the complex conjugation 
and w.r.t. the rotation by
the $r$-th root of unity. 
In particular, it is invariant by $z\mapsto -z$.
Hence, Theorem A(1) implies that
if $R\in {\M}_0$ then $R/\alpha\in {\M}_0$.
If now $g^n(R)=R_0$ for some $R_0\in {\M}_0$, then
$g^{2n}(R/\alpha)=g^n(R)/\alpha=R_0/\alpha\in {\M}_0$.
Let's show that $\partial R\cap J_\infty$
is the base point $\{x_R\}$ of $R$.
First, $\partial R\cap J=\{x_R\}$, from the definition.
If now $x\in \partial R\cap \alpha^n J$,
for some $n$, then
$x/\alpha^n\in \partial (R/\alpha^n)\cap J=\partial G\cap J=
\{x_{G}\}$, for a piece $G$. Hence,
$x=\alpha^n x_G=x_R$.

{\bf (b)} follows from Theorem A(5).
In fact, by the same reason, a stronger statement holds.
Namely, remind that $R_0$ is an 
infinite open ``sector'' 
bounded by the curves $L$ and $-L^\ast$.
Introduce the set 
$\tilde {\M}=\{\tilde R\}$ of all components of
$g^{-n}(R_0)$ and $g^{-n}(R_0^\ast)$
for all $n\ge 1$. Then any two such components
are either disjoint or one covers the other one.
Repeating (a), we get also that if $\tilde R\in \tilde {\M}$ then
$\pm \tilde R/\alpha\in \tilde {\M}$.

Note that any piece $R$ is contained in one and only one 
element $\tilde R\in \tilde {\M}$ such that $R$ and $\tilde R$
have the same boundary ``base points'' 
(remind that if $g^n(R)\in {\M}_0$ then the ``base points'' of $R$
are $x_R$, $x_R^\infty$ s.t. $g^n(x_R)=0$, $g^n(x_R^\infty)=\infty$).

{\bf (c)-(d).} It is enough to prove both statements
replacing $R\in {\M}$ by $\tilde R\in \tilde {\M}$. 

Let's construct a ``machine'' that produces 
new elements of  $\tilde {\M}$.
Consider a component $\tilde R_1$ of $g^{-1}(R_0^\ast)$, namely,
$\tilde R_1=u^\ast(R_0)$. Its boundary consists of the points 
$c, x_0$ and two
arcs $\gamma^\pm$ joining the points $c, x_0$,
such that $g(\gamma^+)=-L$ and $g(\gamma^-)=L^\ast$.
Moreover,
$\gamma^-\cap L=L_0$ (because $g(L_0)=L_\infty^\ast\subset L^\ast$)
while $\gamma^+$ has no common
arcs with $L$. Denote by $\gamma^-_1$ the sub-arc of $\gamma^-$
between the points $c/|\alpha|$ and $x_0$:
$\gamma^-_1=\gamma^-\setminus \overline L_0$ . In turn, 
it is divided into
two arcs: $\beta_1$ and $\beta_2$ such that
$g(\beta_1)=L_0^\ast$ and $g(\beta_2)$ is the arc of $L^\ast$
between $0$ and $c^\ast/|\alpha|$.
Observe that 
\begin{equation}\label{machinearc}
\gamma^+/|\alpha|\cap \gamma^-_1=\beta_1.
\end{equation}
(Proof: $g^2(\gamma^+/|\alpha|)=g(\gamma^+)/\alpha=-L/\alpha=L$;
hence, $g(\gamma^+/|\alpha|)=(\gamma^-)^\ast$ and
$g(\gamma^+/|\alpha|\cap \gamma^-_1)=L_0^\ast=g(\beta_1)$.)

Denote by $\Delta_0$ an open ``triangle'' 
(with the ``vertices'' $x_0, x_0/|\alpha|, c/|\alpha|$) bounded by
the interval $(x_0/|\alpha|, x_0)$, the sub-arc 
$\beta_2=\gamma^-_1\setminus \beta_1$ of $\gamma^-_1$,
and the sub-arc $\gamma^+/|\alpha|\setminus \beta_1$ of $\gamma^+/|\alpha|$.
Denote also by $\Delta$ a (bigger) open ``triangle''
bounded by the interval $(0, x_0)$, the arc $\gamma^-_1$, and
the arc $L\setminus \overline{L_0\cup L_\infty}$
of $L$ between $0$ and $c/|\alpha|$.
It is now easy to check that
$g(\partial \Delta_0)=\partial \Delta^\ast$, i.e.
$\Delta_0=u^\ast(\Delta)$.
Also, $R_2:=\tilde R_1/|\alpha|\subset \Delta$.
We have now a dynamical system (the ``machine'')
consisting of two maps
$A: z\mapsto z/|\alpha|$ and
$B: z\mapsto u^\ast(z)$ of $\Delta$
into itself which produces new elements of $\tilde {\M}$: 
if $\tilde R\in \tilde {\M}$ is a subset
of $\Delta$, then $A(\tilde R)$ and $B(\tilde R)$
are again elements of $\tilde {\M}$.
To study it, note that there is a neighborhood
$U$ of $\overline\Delta$, such that
$B=u^\ast$ extends to a conformal map on $U$,
$\overline{B(U)}\subset U$ and  $\overline{A(U)}\subset U$
(for example, take a Poincare neighborhood of an interval
$(-\e, 1)$ in the slit plane $(\C\setminus \R)\cup (-\e, 1)$).
Then $A, B:\overline \Delta\to \overline \Delta$ are stricts contractions
in the hyperbolic metric of $U$.
Besides, $A(\Delta)$, $B(\Delta)$, and $R_2$ are pairwise disjoint,
and 
$A(\overline\Delta)\cup B(\overline\Delta)\cup \overline R_2=\overline\Delta$.
It follows that if we apply to $R_2$ 
all possible finite compositions
of $A$ and $B$, then the union of the closures of all elements
of $\tilde {\M}$
obtained like this, with the points of the segment
$[0, x_0]$, is just the closure of $\Delta$. 

Let's finish the proof of (c). As we have just seen,
the set $F=\cup \pm(\overline {\Delta\cup \Delta^\ast})$ is filled
by the closed pieces (and the points of the segment $[-x_0, x_0]$).
$F$ also a closed neighborhood of the point
$g^2(0)=1/\alpha\in (-x_0, 0)$.
Pulling $F$ back by a branch of $g^{-2}$ we get
a neighborhood of zero with the same property. 

Let's prove (d). Let $R$ be any piece.
Since $R$ is a preimage of a piece of the depth one,
there indeed exists a minimal $n\ge 0$ such that
$g^n(R)=G/\alpha^k$, for some piece $G$ 
of the depth one and for some $k\ge 0$.
Then $g^n(\tilde R)=\tilde G/\alpha^k$,
for corresponding $\tilde R\supset R$ and $\tilde G\supset G$
(so that $g(\tilde G)$ is either $R_0$ or $R_0^\ast$).
To prove (d2), we find for each (among $2r$)
component $\tilde G$ of $g^{-1}(R_0)$ or $g^{-1}(R_0^\ast)$
a fixed neighborhood $V(\tilde G)$, so that
the branch of $g^{-n}$ from $\tilde G/\alpha^k$ onto $\tilde R$
extends to a map onto $V(\tilde G)/\alpha^k$
 
If $\tilde G$ does not touch the real axis, such a neighborhood
obviously exists.
The remaining possibility is that $\tilde G$ is one of
$\pm \tilde R_1$, $\pm \tilde R_1^\ast$ (see the definition of 
$\tilde R_1$ above).

For example, let $\tilde G=\tilde R_1$.
Let us fix a neighborhood $V$ of the closure of
$\tilde R_1$, so that $V$ intersects 
the real axis along, say, the interval $(g^3(0), 1+\epsilon)$,
with $\epsilon>0$ fixed and small enough.
Assume that
the branch of $g^{-n}$ from $\tilde R_1/\alpha^k$ onto $\tilde R$
does not extend to a map on $V/\alpha^k$,
i.e., $g^{-n}=g^{-(n-m)}\circ g^{-m}$
where $g^{-m}(V/\alpha^k)$ covers $0$ for the first time. Then, 
as the real Feigenbaum dynamics shows us, $m=2^k$, and, moreover, 
the branch $g^{-m}(z)=g^{-2^k}(z)=(1/\alpha^k)g^{-1}(\alpha^k z)$,
hence,
the element $g^{n-m}(\tilde R)=g^{-m}(\tilde R_1/\alpha^k)$
has the form $\tilde G_1/\alpha^{k+1}$, for a component $\tilde G_1$
of either $g^{-1}(R_0)$ or $g^{-1}(R_0^\ast)$,
a contradiction with the minimality of $n$.
 
Thus (d1)-(d2) are proved.

To complete the proof show that if
$R_{m}\subset R_{m-1}\subset...\subset R_1$ are different pieces,
then $diam(R_m)\le C\lambda^m diam(R_1)$, for universal
$C>0, \lambda<1$.
First, observe, that any piece $G$ of depth one contains at least
two (in fact, infinite number) of disjoint pieces of level $2$
(i.e. which are not contained in any other piece except $G$).
Hence, by (d1)-(d2), there exists $\rho<1$, so that
$area(R')/area(R)\le \rho$ for arbitrary different pieces $R'\subset R$.
It follows, $area(R_m)\le \rho^m area(R_1)$.
Second, again by (d1)-(d2) and by compactness of an appropriate
family of univalent maps, there is $C$, such that
$C^{-1}\le (diam(R))^2/area(R)\le C$ for all pieces $R$.
Therefore, $diam(R_m)\le C\rho^{m/2} diam(R_1)$.  

\

The proof of (c)-(d) yields the following 
\begin{com}\label{levelone}
We have seen that the ``machine'' produces
the covering of the triangle $\Delta$
by the pairwise disjoint pieces, and each of them
touches the interval $(0, x_0)$.
It follows that every piece $R$ of the level $1$
touches at its ``base point'' $x_R$
one of the lines ${\R}\exp(i\pi j/r)$, $j=0,1,...,r-1$.
\end{com}
\begin{com}\label{l(x,y)}
Since every point of $\overline\Delta\cap \R^+$
is a limit of real ``base points'' $x_{\tilde R}$
of elements $\tilde R$ produced by the ``machine'',
for every point $x\in \overline\Delta$, which
lies either in ${\R}^+$ or in the boundary
of a piece of level one, there exists a simple curve
$l(x,c/|\alpha|)$ which joins the points
$x$ and $c/|\alpha|$, lies in the (open)
``triangle'' $\Delta$ and such that $l(x,c/|\alpha|)$
consists of closed sub-arcs of boundary curves of
pieces of level one.
\end{com}
\begin{com}\label{levelzero}
The "ray" $g^2(L\cap g^{-2}(\Omega))$, in a neighborhood
of the landing point $g^2(0)$, consists of (closed,
i.e. with their end points) arcs of boundaries
of pieces of level one (follows from Theorem A(5)). 
This sequence of arcs accumulates only
at the point $g^2(0)$, as it follows from Theorem B (d)
and the previous Comment.
Pulling this picture back, we see 
the same is true for the "ray" $L$ (with
the obvious replacement of $g^2(0)$ by $0$): 
$L$ consists of a sequence of (closed) arcs of
$\partial R$, where $R$ are some pieces of level one,
which touch the straight ray
${\R}^+\exp(i\pi/r)$.
Remind also that at the same time $L$ consists of
a sequence of (closed) arcs of $\partial R$ with ends at the points
$c/|\alpha|^n, n\in {\Z}$, where $R$ are some pieces of level one,
which touch the real ray ${\R}^+$.
\end{com}

\section{Some conclusions}
Let us summarize some useful facts that follow
from the previous section.
\subsection{The structure of a piece}
Let $R\in {\M}^{ext}$ be any piece. 
By the definition, $g^m(R/\alpha^j)=G$,
for some $m, j\ge 0$, where $G$ is a piece of the depth zero, i.e.
one of the infinite ``sectors''
$R_{0,1}, R_{0,2},..., R_{0,r-1}$,
$R_{0,1}^\ast, R_{0,2}^\ast,..., R_{0,r-1}^\ast$
(see the definition before Theorem B).
The $G$ contains one and only one straight ray of the form
${\R}^{+}\exp(i\pi j/r), j=\pm 1,...,\pm(r-1)$.
Denote it by $\gamma_G$.
The following facts have been proved.

{\bf 1.} The boundary of $R$, $\partial R=L(R)\cup \{x_R\}$,
where $x_R\in B$ is a ``base point'', i.e. $g^m(x_R/\alpha^j)=0$,
and $L(R)$ is a simple (open) curve 
which consists of arcs of preimages of the curve $L$.
Hence, $L(R)\subset {\C}\setminus J_\infty$.

{\bf 2.} Define the ``vein'' of $R$, $\gamma_R$, so that
$g^m(\gamma_R/\alpha^j)=\gamma_G$. 
Then $\gamma_R$ is an analytic arc joining
the ``base points'' $x_R, x_R^\infty$
of $R$. Moreover,
$\gamma_R\subset J_\infty$, and the length of $\gamma_R$ is finite
(Comment~\ref{lens} and Theorem B (d)).

{\bf 3.} Let $n\ge 1$ be the level of
the piece $R$.
Let us look at the set ${\M}^{ext}(R)$ of all pieces of level $n+1$
which are contained in $R$. Then:

\begin{itemize}
\item
(3a) they form a partition of $R$, i.e. 
they are pairwise disjoint, and the union
of their closures $\overline P$ ($P\in {\M}^{ext}(R)$)
united with the (closed)
``vein'' $\gamma_R\cup \{x_R\}\cup \{x_R^\infty\}$
is the closure of $R$,
\item
(3b) each $P\in {\M}^{ext}(R)$ touches
the ``vein'' $\gamma_R$ at the point $x_P$ (Comment~\ref{levelone});
the boundaries of $P$, $R$ are either
disjoint or have a common non-trivial arc,
\item
(3c)  the latter is true for any two
pieces $P, P'\in {\M}^{ext}(R)$, with an
extra possibility that $P, P'$ can lie on different sides 
of the ``vein'' and have the same ``base point'' $x_P=x_{P'}$, 
\item
(3d) all $P\in  {\M}^{ext}(R)$ as well as $R$
have roughly the same shape (see Theorem B (d)).
\item
(3e) the "base points" $x_R, x_R^\infty$ divide
$\partial R$ into two simple curves: $L^{\pm}(R)$
(in other words, $L^+(R)\cup L^-(R)\cup \{x_R^\infty\}=
L(R)$). As it follows from Comment~\ref{levelzero},
either curve $L^{\pm}(R)$ consists of a sequence
$\{L^{\pm}_n(R)\}_{n=-\infty}^{+\infty}$ of (closed) sub-arcs
of boundaries of pieces $R_n^\pm$ from ${\M}^{ext}(R)$
(so that $x_{R_n^+}=x_{R_n^-}$),
which accumulate at the points $x_R, x_R^\infty$:
$\lim_{n\to -\infty} L^{\pm}_n(R)=\{x_R\}$,
$\lim_{n\to +\infty} L^{\pm}_n(R)=\{x_R^\infty\}$.

For every large enough integer $m$, denote by
$W_m(R)$ the domain bounded by the arcs $L^\pm_n, n\ge m$,
the points $x_R^\infty$ and $x_{R_m^+}=x_{R_m^-}$,
and arcs of $\partial R_m^\pm$ from $x_{R_m^\pm}$ to $\partial R$.
\item
(3f) if $y$ is an end point of an arc $L^\pm_n(R)$,
then, for some non-negative $n, j$,
the point $g^n(y/\alpha^j)$
belongs to the set $\{ c\exp(i2\pi j/r), c^\ast\exp(i2\pi j/r),
j=0,...,r-1\}$.
Indeed, it is enough to prove it when $R$ is of
depth zero. Note that $y$ belongs to the
boundaries of two different pieces.
Then iterating it further we apply 
Theorem A (5) and get the statement. 
\end{itemize}
\subsection{Neighborhoods of points}\label{W}
Fix a point $z\in \partial M^{ext}\setminus Y$.
That is, for some non-negative $n, j$, the point
$w=w(z)=g^n(z/\alpha^j)$ is either infinity or belongs to 
one of the curves
$L\exp(i 2\pi j/r), L^\ast\exp(i 2\pi j/r)$ 
($0\le j\le r-1$). We want to describe a basis
of neighborhoods at $z$ as union of pieces.
Consider two cases.
\begin{itemize}
\item
$z$ is not a ``base point'' $x_R^\infty$ of any
piece $R$. Then, by Comment 4, the corresponding $w$ is a boundary
point of at least $2$ and at most $3$ closed pieces of level $1$,
such that their union forms a closed neighborhood of $w$.
Returning to $z$ and using (3e), we find
a sequence of closed neighborhoods $W_n(z)$ of $z$
shrinking to $z$, such that
each $W_n(z)$ is the union of at least $2$
and at most $3$ closed pieces.
\item
$z$ is a ``base point'' $x_R^\infty$
of a piece $R$. It follows from the definition
of $x_R^\infty$, that then there are exactly $r-1$ pieces
$R_j, j=1,...,r-1$, so that $z=x_{R_j}^\infty$.
Moreover, all $R_j$ are of the same depth
(and, in fact, of the same level).
Similarly to the first case, we then get 
a sequence of closed neighborhoods $W_n(z)$ of $z$
shrinking to $z$, such that
each $W_n(z)$ is the union of at least $1$
and at most $2$ closed pieces AND $r-1$ closures
of domains $W_n(R_j), j=1,...,r-1$ defined in (3e)
(see the previous subsection).
\end{itemize}
\section{Theorem C and Corollary}\label{theoc}
\subsection{Proof of Theorem C}
If $z\in J^{ess}$, {\bf (a)} follows directly from Theorem B(b)-(c).
If $z\in \partial M^{ext}\setminus Y$, then {\bf (a)}
follows from Section~\ref{W}.

To prove {\bf (c)} we need 

{\bf Fact.} {\it There is $C'$ such that, for any
piece $R\in {\M}$ and any point 
$w\in \overline{R}$, there exists
a rectifiable curve $\sigma(w, x_R)\subset R\cap J_\infty$
(in fact, consisting of closed (arcs of) ``veins''
of different pieces), which joins $w$ and the ``base point'' $x_R$
and of length less than $C'$diam$R$.}

Indeed, let $R$ be of level $n$; if $z\notin Y$, then, by {\bf (a)},
$\{w\}=\cap_{k\ge n} \overline{R_k(w)}$
where $R_k(w)$ is a piece of level $k$
s.t. $w\in \overline{R_k(w)}$.
Moreover, for any $k\ge n$, $R_{k+1}(w)$ is contained
in $R_{k}(w)$ and touches
the ``vein'' of $R_k(w)$. With this representation,
the curve $\sigma(w, x_R)$ arises naturally:
we start at $x_R$ and follow the ``vein'' $\gamma_R$ until
the base point $x_{R_{n+1}(w)}$, then follow the ``vein''
$\gamma_{R_{n+1}(w)}$ until the base point $x_{R_{n+2}(w)}$
etc. 
Now Fact is a corollary of Theorem B (d).
If $z\in Y$, then $y$ belongs to a ``vein''
of a piece $R_m$ contained in $R_n(w)$.
Then we repeat the above construction going
along ``veins'' 
from the piece $R_m$ until the base point 
$x_R$.

Then {\bf (c)} follows from the Fact and Theorem B (c).

It remains to prove {\bf (b)} and {\bf (b')}.
The {\bf existence} of the curve $l(x, y)$ ($x, y\in \overline \Pi$)
follows from (i)-(vii):
\begin{itemize}
\item
(i) if $l(x, y)$ is constructed then one can
set $l(|\alpha|^j x, |\alpha|^j y)=|\alpha|^j l(x, y)$,
$j\in {\Z}$.
\item
(ii) it is enough to construct $l(x, c)$, for any
$x\in \overline \Pi$, because then $l(x, y)$ 
will be a composition of the curves $l(x, c), l(c, y)$
with cancelling out arcs passing in opposite directions. 
\item
(iii) if $x\in \R^+$ or $x\in \partial R$,
where $R$ is a piece of level one,
which lies in the closed ``sector'' $S$ bounded
by the curves $L$ and ${\R}^+$, then 
using (i) and Comment~\ref{l(x,y)}, we see
that the curve $l(x, c)$ exists,
belong to $S$ and it goes inside some of the curves $L(R)$
where $R$ runs over the pieces of level one.
\item
(iv) let $x\in \overline G$, for a piece $G$ of
{\bf depth} zero (since $x\in \overline \Pi$,
it means that $x$ belongs to the closed
``sector'' bounded by $L$ and $\exp(2\pi i/r)\R^+$).
Moreover, assume that $x$ lies either on
the ray $\exp(2\pi i/r)\R^+$ or in $\partial R$,
for a piece $R$ of level one.
Consider the point $g(x)$. It belongs to the ``sector'' $S$
(defined in (iii)).
Using (i) and (iii), we obtain the curve
$l(g(x), c)$, and then define the curve $l(x, g^{-1}(c))=
g^{-1}(l(g(x), c))$, with a branch
$g^{-1}: g(G)\to G$. Extending the curve $l(x, g^{-1}(c))$
along the boundary of $G$ to infinity, we get
a curve $l(x, \infty)$ which lies in $\overline G$
and consists of arcs of the curves $L(R)$
where $R$ are pieces of level one.

Cases (iii)-(iv) cover the points $x\in \overline \Pi$ belonging either to
two rays $\partial \Pi$ or to the boundary of a piece
of level one.

\item
(v) now we can go into deeper levels.
Let $x\in \overline R_n$, for some $R_n\in {\M}^{ext}$,
a piece of level $n$, and such that $x$ lies either
on the ``vein'' $\gamma_{R_n}$ or
in the boundary $\partial R_{n+1}$
of a piece $R_{n+1}$ of level $n+1$.
Using (i) one can assume that $R_n\in {\M}$.
If now $m\ge 0$ is so that
$g^m(R_n)=G\in {\M}_0$, then we can apply
(iv) to the point $g^m(x)\in \overline G$.
Going back to $\overline R_n$ by a branch of 
$g^{-m}$, we get the curve $l(x, x_{R_n}^\infty)\subset
\overline R_n$. It consists of arcs of 
curves $L(R)$, for a collection of pieces
$R\subset R_n$ of level $n+1$.
\item
(vi) let $x\in Y\setminus \partial \Pi$. Then $x$ lies in a ``vein''
of a closed piece $R_n$ of a maximal {\it level} $n$. Then use (v)
to get $l(x, x_{R_n}^\infty)$. Since $x_{R_n}^\infty$
lies in a closed piece $R_{n-1}$ of level $n-1$, we get
$l(x_{R_n}^\infty, x_{R_{n-1}}^\infty)$ which
composing with $l(x, x_{R_n}^\infty)$ gives
$l(x, x_{R_{n-1}}^\infty)$ etc. until we get
$l(x, x_R^\infty)$, where $x_R^\infty$
is a base point of a piece $R$
of level one. Then apply (iv).
Above considerations treat also the case when
$x\in \partial {\M}^{ext}$.
\item
(vii) finally, let $x\in J^{ess}$.
Then $\{x\}=\cap_{n=1}^\infty R_n(x)$,
and $l(x, c)$ can be defined as the
composition of curves $l(x_{R_n(x)}^\infty, x_{R_{n+1}(x)}^\infty)$
for all $n$.
\end{itemize}

Let's prove the uniqueness. If $x,y$ are not in $J_\infty$, then
$l(x, y)$ is unique because otherwise
(as in the proof of Corollary~\ref{anal}) 
$J_\infty$ would not be connected. Let now $x\in |\alpha|^n J\setminus Y$.
If there are two different simple curves joining
$x, y$ and satisfying (b1)-(b2) of Theorem C, then
again because $J_\infty$ is connected, for some
$m\ge n$, the set $|\alpha|^m J\setminus \{x\}$ has at least
two components. This is equivalent to say that
$x\in |\alpha|^m Y_0$, a contradiction.

{\bf (b')}: (i) is clear from
the construction, (ii): notice that, by the
construction of the curves $l(x, y)$,
any two such curve are either disjoint, or
intersected along
an arc, so that each its end point belongs to
two different pieces. It remains to apply
the property (3f) of Sect. 4.1. 

\subsection{Proof of Corollary 1.1 (b2)}
We need to show that the curves $l(x, y)\to z$
whenever $x, y, z\in \partial W$ and $x, y\to z$.

If $z\in J^{ess}$, then $z$ is the intersection of nested
sequence of (open) pieces $R_n(z)$, so that
$R_n(z)$ is a neighborhood of $z$. 
Then the statement follows from Theorem C (a).

If $z\in \partial M^{ext}\setminus Y$, then 
we use the sequence of closed neighborhoods $W_n(z)$ of $z$
constructed in Section~\ref{W}:
$l(x, y)\subset W_n(z)$ whenever $x, y\in  W_n(z)$.
Since $W_n(z)$ shrink to $z$, we get the result.
\section{Appendix: proof of Lemma~\ref{meas} and Corollary~\ref{frost}}
Remind that $I$ is the limit set of the system
of maps $S=\{\ph_i:X\to X\}_{1\le i\le 3}$.

{\bf Proof of Lemma~\ref{meas}.} We need to find 
$0<h\le 2$ and a non-atomic probability 
measure $\mu$ with its support $supp(\mu)=I$,
which is $h$-conformal, i.e., 
$\mu(\ph_i(A))=\int_A |\ph_i'(x)|^h d\mu(x)$ for
$i=1,2,3$ and any Borel set $A\subset X$. 
There are two ways of proving such statements.
One method goes back to Bowen and Ruelle, and is developed e.g.
in~\cite{MU}; first, it is shown the existence
of so-called semi-conformal measure,
which is a fixed point of the corresponding adjoint
Ruelle transfer operator, then the condition (5) of Lemma~\ref{ifs}
can be used to prove that this semi-conformal measure
is actually conformal (see the proof of
Lemma 3.10 of~\cite{MU}).
Here we follow the second, more direct approach of~\cite{sum}, 
where the existence 
of conformal measure for rational maps is proved.
Main difference with~\cite{sum} is that
the ``expanding'' system of inverse maps $\ph_i^{-1}$
does not glue to a global map. On the other hand, the set
of points where it does not happen, 
$\ph_i(X)\cap \ph_j(X), i\not=j$, is ``small''
because of Lemma~\ref{ifs} (5).

Fix a point $z\in Int(X)\setminus I$.
For every $s>0$, define the series
$$P(s):=\sum_{\e\in \Delta} |\ph_\e'(z)|^s,$$
where $\Delta$ is the set of all finite words $\e$.
Considering the area case ($s=2$), we see that $P(s)$ is 
finite for $s\ge 2$. Define
$$h=\inf\{s: P(s)<\infty\}.$$ 
Then $h>0$ because $\inf \{i, x\in X: |\ph_i'(x)|\}>0$.
Pick up a positive function $\tau(t), t>0$, such that:

(i) $\tau(t)\to +\infty$ as $t\to 0$ in such a way, that
for all $\delta>0$, $0<\lambda< \infty$,
there is $t_0(\delta, \lambda)$ so that
$\tau(\lambda t)/\tau(t)\in (1-\delta, 1+\delta)$
for $0<t\le t_0$,

(ii) $P_\tau(h):=
\sum_{\e\in \Delta} \tau(|\ph_\e'(z)|) |\ph_\e'(z)|^h=\infty$

(If $P(h)=\infty$, there is no need for introducing the function $\tau$).

Notice that (i) implies that $\tau(t)<t^{-a}$ for every $a>0$
and all $t$ small enough. Therefore, 
$P_\tau(s):=\sum_{\e\in \Delta} \tau(|\ph_\e'(z)|) |\ph_\e'(z)|^s<\infty$
for $s>h$.

For every $s>h$, define a probability measure
\begin{equation}\label{mus}
\mu_s={1\over P_\tau(s)}
\sum_{\e\in \Delta} \tau(|\ph_\e'(z)|) |\ph_\e'(z)|^s \delta_{\ph_\e (z)}.
\end{equation}
Let $\mu$ be any weak limit of $\mu_s$ as $s$ desreases to $h$.
Then (ii) implies that $\mu$ is supported on $I$.
Fix $i\in \{1,2,3\}$.
First consider any point $x_0\in I\cap Int(X)$ and 
consider a small enough neighborhood $U$ of $x_0$.
Then $\ph_j(x)\in \ph_i(U)$, for some $x\in X$, implies 
that $j=i$ and $x\in U$.
Because of (i) and since the length of words $\e$, such that
$\ph_\e(z)\in U$ , tends to infinity as $U$ shrinks to $x_0$,
we derive that for any $s>h$,
\begin{equation}\label{pr1}
\mu_s(\ph_i(U))/|\ph_i'(x_0)|^s\mu_s(U)\in [K^{-1}, K],
\end{equation}
where the constant $K>1$ tends to $1$ 
uniformly in $s$ as $U$ shrinks to $x_0$.

Let's show that any $x\in I$ is not an atom of $\mu$.
Assume the contrary.
If $x=\cap_{n=1}^\infty \ph_{\e|n}(Int(X))$
for an infinite word $\e$, then applying~(\ref{pr1})
finitely many times to the points
$\ph_{\e|n}^{-1}(x)$ we come to a contradiction
with the fact that the measure $\mu$ is finite.
As for the rest of the points $x$, the same argument shows
that one can assume from
the beginning that $x\in I\cap \cup_{1\le j\le 3}\partial \ph_j(X)=
\{a_j, 1\le j\le 4\}$.
If $x$ is either $a_1$ or $a_4$, then $x$
is the attracting fixed point of either $\ph_1$ or $\ph_3$.
If, say, $x=a_1$, then again, for small neighborhoods
$U$ of $a_1$,
$\ph_j(y)\in \ph_1(U)$, for some $y\in X$, implies 
that $j=1$ and $y\in U$. Then we can apply the same argument as above
to the map $\ph_1$, and get a contradiction.
Consider the remaining cases when $x$ is either $a_2$ or $a_3$,
say, $x=a_2$. Since $a_2=\ph_3(a_1)=\ph_1(a_4)$,
for a small neighborhood
$U$ of $a_2$,
$\ph_j(y)\in U$, for some $y\in X$, implies 
that $j$ is either $3$ or $1$, and $y$ lies in a neighborhood
of $a_1$, $a_4$ resp.  It shows that either $a_1$ or $a_4$
is an atom, a contradiction.

Since $\mu$ has no atoms and $I\cap \cup_{1\le i\le 3} \ph_i(\partial X)$
is finite, it is easy to see that
the measure $\mu^\ast$ defined by
$\mu^\ast(A)=\mu(\ph_i(A))$, for any Borel $A\subset X$,
is absolutely continues
w.r.t. $\mu$. By the same reason, $x_0$ as in~(\ref{pr1}) is a typical
point of $\mu$. Passing to the limit 
in~(\ref{pr1}) as $s\to h$ and $U$ shrinks to $x_0$ 
we have $d\mu^\ast/d\mu=|\ph_i'|^h$ $\mu$-a.e.

{\bf Proof of Corollary~\ref{frost}.} Prove that
$h=HD(I)$ and the $h$-Hausdorff measure of $I$ is finite
and positive. This will follow in a usual way (see 
e.g.~\cite{MU}, Lemma 3.14)
from the finiteness of the system $S$ and the existence
of the conformal measure. 
By a general fact from geometric measure theory,
see e.g.~\cite{ma}, it is enough to show that there exists a constant $C>0$,
such that, for every $x\in I$ and every $0<r<1$,
\begin{equation}\label{frostin}
C^{-1}\le {\mu(B(x,r))\over r^h}\le C.
\end{equation}
Notice that by the bounded distortion property
there is $C_0$ such that, for every finite word $\tau$,
$|\tau|=n$, for every $y\in X$, we have
$C_0^{-1}\le diam(\ph_\tau(X))/|\ph_\tau'(y)|\le C_0$,
$C_0^{-1}\le area(\ph_\tau(X))/ diam(\ph_\tau(X))^2\le C_0$,
$C_0^{-1}\le diam(\ph_\tau(X))/diam(\ph_{\tau|n-1}(X))\le C_0$.
and $C_0^{-1}\le \mu(\ph_\tau(X))/ diam(\ph_\tau(X))^h\le C_0$.
Let $x=\cap_{n=1}^\infty \ph_{\e|n}(X)$.
Let $m$ be minimal such that $\ph_{\e|m}(X)\subset B(x,r)$.
From the minimality of $m$ and because
$x\in \ph_{\e|m}(X)$, $diam(\ph_{\e|m}(X))\ge r/C_0$, hence,
$\mu(B(x,r))\ge diam(\ph_{\e|m}(X))^h/C_0\ge r^h/C_0^{h+1}$.
To prove the opposite inequality, 
consider a collection $\Sigma$
of sets of the form $\ph_w(X)$, such that
$B(x,r)\cap \ph_w(X)\not=\emptyset$, $\ph_w(X)\subset B(x,2r)$, and
the cardinalities $|w|$
are minimal (i.e. if $\ph_w(X)\in \Sigma$, then
$\ph_{w||w|-1}(X)$ is not a subset of $B(x, 2r)$). 
Clearly, these sets form a cover of $B(x,r)\cap I$,
and their interiors are pairwise disjoint. 
From the minimality of $|w|$, we conclude
that $diam(\ph_w(X))\ge r/C_0$. Hence,
$4\pi r^2=area B(x,2r)\ge \sum_{A\in \Sigma} area(A)\ge
C_0^{-1}\sum_{A\in \Sigma} diam(A)^2\ge N r^2/C_0^3$,
where $N$ is the number of elements in $\Sigma$. Therefore,
$N\le 4\pi C_0^3$, and, finally,
$\mu(B(x,r))\le \sum_{A\in \Sigma} \mu(A)\le
C_0\sum_{A\in \Sigma} diam(A)^h\le N C_0 (4r)^h\le
4^{1+h}\pi C_0^4 r^h$.

\

{\bf Acknowledgments.} I would like to thank
Feliks Przytycki, Juan Rivera-Letelier,
Sebastian van Strien, and Greg Swiatek
for discussions, and the referee for
very useful remarks that led to a revision
of the paper.

\end{document}